\input amstex
\documentstyle{amsppt}
\magnification=\magstep1                        
\hsize6.5truein\vsize8.9truein                  
\NoRunningHeads
\loadeusm

\document

\topmatter

\title
Tur\'an-type reverse Markov inequalities for polynomials with restricted zeros 
\endtitle

\rightheadtext{Tur\'an-type reverse Markov inequalities for polynomials with restricted zeros}

\author 
Tam\'as Erd\'elyi 
\endauthor

\address Department of Mathematics, Texas A\&M University,
College Station, Texas 77843, College Station, Texas 77843 \endaddress

\endthanks


\date September 21, 2019 
\enddate

\email terdelyi\@math.tamu.edu
\endemail

\abstract
Let ${\Cal P}_n^c$ denote the set of all algebraic polynomials of degree at most $n$ with complex coefficients.
Let
$$D^+ := \{z \in {\Bbb C}: |z| \leq 1, \enskip \text {\rm Im}(z) \geq 0\}\,.$$
For integers $0 \leq k \leq n$ let ${\Cal F}_{n,k}^c$ be the set of all polynomials $P \in {\Cal P}_n^c$ having 
at least $n-k$ zeros in $D^+$. Let 
$$\|f\|_A := \sup_{z \in A}{|f(z)|}$$
for complex-valued functions defined on $A \subset {\Bbb C}$.
We prove that there are absolute constants $c_1 > 0$ and $c_2 > 0$ such that
$$c_1 \left(\frac{n}{k+1}\right)^{1/2} \leq 
\inf_{P}{\frac{\|P^{\prime}\|_{[-1,1]}}{\|P\|_{[-1,1]}}} \leq c_2 \left(\frac{n}{k+1}\right)^{1/2}$$
for all integers $0 \leq k \leq n$, where the infimum is taken for all $0 \not\equiv P \in {\Cal F}_{n,k}^c$ 
having at least one zero in $[-1,1]$. This is an essentially sharp reverse Markov-type inequality 
for the classes ${\Cal F}_{n,k}^c$ extending earlier results of Tur\'an and Komarov from the case $k=0$ 
to the cases $0 \leq k \leq n$.  
\endabstract

\endtopmatter

\head 1. Introduction and Notation \endhead
Let ${\Cal P}_n$ denote the set of all algebraic polynomials of degree at most $n$ with real coefficients
Let ${\Cal P}_n^c$ denote the set of all algebraic polynomials of degree at most $n$ with complex coefficients. 
Let 
$$\|f\|_A := \sup_{z \in A}{|f(z)|}$$
for complex-valued functions defined on $A \subset {\Bbb C}$.
In 1939 Tur\'an [T-39] proved that 
$$\|P^\prime\|_{[-1,1]} \geq \frac{\sqrt n}{6} \, \|P\|_{[-1,1]} \tag 1.1$$
for all $P \in {\Cal P}_n^c$ of degree $n$ having all their zeros in the interval $[-1,1]$. The examples 
$P(x) = (x^2-1)^m$ and $P(x) = (x^2-1)^m(x+1)$ show that Tur\'an's reverse Markov-type inequality (1.1) is 
essentially sharp, even  though the multiplicative constant $1/6$ in (1.1) is not the best possible. 
Note that the best possible multiplicative constant $c = c_n$ in (1.1) had been found by Er\H od [E-39]. 
Another simple observation of Tur\'an [T-39] is the inequality
$$\|P^\prime\|_D \geq \frac n2 \, \|P\|_D \tag 1.2$$
for all $P \in {\Cal P}_n^c$ of degree $n$ having all their zeros in the closed unit disk $D \subset {\Bbb C}$.
Malik [M-69], Govil [G-73], and Govil and Mohapatra [GM-99, Section 4] established extensions of (1.2) proving that 
$$\|P^\prime\|_D \geq \frac{n}{1+R} \, \|P\|_D$$  
for all $P \in {\Cal P}_n^c$ of degree $n$ having all their zeros in the disk $D(0,R) \subset {\Bbb C}$
of radius $R \leq 1$ centered at $0$,  and
$$\|P^\prime\|_D \geq \frac{n}{1+R^n} \, \|P\|_D$$
for all $P \in {\Cal P}_n^c$ of degree $n$ having all its zeros in the disk $D(0,R) \subset {\Bbb C}$
of radius $R \leq 1$ centered at $0$.

Let $\varepsilon \in [0,1]$ and let $D_\varepsilon$ be the ellipse
of the complex plane with large axis $[-1,1]$ and small axis
$[-i\varepsilon,i\varepsilon]$. Let ${\Cal P}_n^c(D_\varepsilon)$ denote
the collection of all $P \in {\Cal P}_n^c$ of degree $n$ having all their zeros in $D_\varepsilon$.
Extending Tur\'an's reverse Markov-type inequality (1.1), Er\H od [E-39, III. t\'etel] proved that
$$c_1(n\varepsilon + \sqrt{n}) \leq
\inf_{P}{\frac{\|P^\prime\|_{D_\varepsilon}}{\|P\|_{D_\varepsilon}}} \leq c_2(n\varepsilon + \sqrt{n})\,,$$
where the infimum is taken for all $P \in {\Cal P}_n^c(D_\varepsilon)$. Levenberg and Poletcky [LP-02] 
rediscovered this beautiful result. In [LP-02] they also proved that   
$$\frac{\sqrt{n}}{20 \, \text {\rm diam} K} \leq \inf_{P}{\frac{\|P^\prime\|_K}{\|P\|_K}}$$
for all compact convex set $K \subset {\Bbb C}$, where the infimum is taken for all $P \in {\Cal P}_n^c$ 
of degree $n$ having all their zeros in $K$.

Let $\varepsilon \in [0,1]$ and let $S_\varepsilon$ be the diamond of the complex plane with diagonals $[-1,1]$ and
$[-i\varepsilon,i\varepsilon]$. Let ${\Cal P}_n^c(S_\varepsilon)$ denote
the collection of all $P \in {\Cal P}_n^c$ of degree $n$ having all their zeros in $S_\varepsilon$.
It has been proved in [E-07] that there are absolute constants 
$c_1>0$ and $c_2>0$ such that
$$c_1(n\varepsilon + \sqrt{n}) \leq
\inf_{P}{\frac{\|P^\prime\|_{S_\varepsilon}}{\|P\|_{S_\varepsilon}}} \leq
c_2(n\varepsilon + \sqrt{n})\,,$$
where the infimum is taken for all $P \in {\Cal P}_n^c(S_\varepsilon)$ with the property
$$|P(z)| = |P(-z)|\,, \qquad z \in {\Bbb C}\,,$$
or where the infimum is taken for all $P \in {\Cal P}_n^c(S_\varepsilon)$ 
with real coefficients. It is an interesting question whether or not the lower bound 
in the above inequality holds for all $P \in {\Cal P}_n^c(S_\varepsilon)$.
Another result in [E-07] shows that this is the case at least when $\varepsilon = 1$, 
that is, there are absolute constants $c_1>0$ and $c_2>0$ such that
$$c_1n \leq \inf_{P}{\frac{\|P^\prime\|_{S_1}}{\|P\|_{S_1}}} \leq c_2n\,,$$
where the infimum is taken for all (complex) $P \in {\Cal P}_n^c(S_1)$.
Motivated by the above results R\'ev\'esz [R-06b] established the right order Tur\'an-type 
reverse Markov inequalities on convex domains of the complex plane. 
His main theorem contains the above mentioned results in [E-07] as special cases. It states that
$$\frac{\|P^\prime\|_K}{\|P\|_K} \geq c(K)n \qquad \text {\rm with} \qquad c(K) = 0.0003 \frac{w(K)}{d(K)^2}\,,$$
for all $P \in {\Cal P}_n^c$ of degree $n$ having all their zeros in a bounded convex set $K \subset {\Bbb C}$,  
where $d(K)$ is the diameter of $K$ and
$$w(K) := \min_{\gamma \in [-\pi,\pi]}{\left( \max_{z \in K}{\text{\rm Re}(ze^{-i\gamma})} -
\min_{z \in K}{\text{\rm Re}(ze^{-i\gamma})} \right)}$$
is the minimal width of $K$. R\'ev\'esz' proof is also elementary, but rather subtle.
Further reverse Markov-type inequalities may be found in the recent papers 
[R-06a] and [R-13] by R\'ev\'esz and [GR-17a] and [GR-17b] by Glazyrina and R\'ev\'esz. 

G.G. Lorentz, M. von Golitschek, and Y. Makovoz devotes Chapter 3 of their book [LG-96] to incomplete polynomials.
E.B. Saff and R.S. Varga were among the researchers having contributed significantly to this topic.
See [BCS-78], [SV-79], and [SV-81], for instance.   
Reverse Markov- and Bernstein type inequalities were first studied by P. Tur\'an [T-39] and J. Er\H od [E-39]
in 1939 (see also [E-06]). The research on Tur\'an and Er\H od type reverses of Markov- and Bernstein-type inequalities
got a new impulse suddenly in 2006 in large part by the work of Sz. R\'ev\'esz [R-06b], and
several results have been published on such inequalities in recent years, see [E-07], [E-09], [EH-15], [GR-17a], [GR-17b], 
[LP-02], [k-04], [MM-94], [R-06a], [R-13], [XZ-99], and [Z-95], for example.

Let ${\Cal P}_{n,k}$ be the set of all algebraic polynomials, with real coefficients, of degree at most $n+k$ 
having at least $n+1$ zeros at $0$. That is, every $P \in {\Cal P}_{n,k}$ is of the form
$$P(x) = x^{n+1}R(x)\,, \qquad R \in {\Cal P}_{k-1}\,.$$
Let
$$V_a^b(f) := \int_a^b{|f^{\prime}(x)| \, dx}$$
denote the total variation of a continuously differentiable function $f$ on an interval $[a,b]$.
In [E-19] a question [EI-18] asked by A. Eskenazis and P. Ivanisvili related to their paper [EI-19] is 
answered by proving that there are absolute constants $c_3 > 0$ and $c_4 > 0$ such that
$$c_3 \frac nk \leq \min_{0 \not\equiv P \in {\Cal P}_{n,k}}{\frac{\|P^{\prime}\|_{[0,1]}}{V_0^1(P)}}  
\leq \min_{0 \not\equiv P \in {\Cal P}_{n,k}}{\frac{\|P^{\prime}\|_{[0,1]}}{|P(1)|}} \leq c_4 \left( \frac nk +1 \right)$$
for all integers $n \geq 1$ and $k \geq 1$. Here $c_3 = 1/12$ is a suitable choice.

In [E-19] we also proved that there are absolute constants $c_3 > 0$ and $c_4 > 0$ such that
$$\split c_3 \left(\frac nk\right)^{1/2} 
& \leq \min_{0 \not\equiv P \in {\Cal P}_{n,k}}{\frac{\|P^{\prime}(x)\sqrt{1-x^2}\|_{[0,1]}}{V_0^1(P)}} \cr  
& \leq \min_{0 \not\equiv P \in {\Cal P}_{n,k}}{\frac{\|P^{\prime}(x)\sqrt{1-x^2}\|_{[0,1]}}{|P(1)|}} 
\leq c_4 \left(\frac nk + 1\right)^{1/2} \cr \endsplit$$
for all integers $n \geq 1$ and $k \geq 1$. Here $c_4 = 1/8$ is a suitable choice.

Let 
$$D^+ := \{z \in {\Bbb C}: |z| \leq 1, \enskip \text {\rm Im}(z) \geq 0\}\,.$$
In [K-19] Komarov proved that 
$$\|P^\prime\|_{[-1,1]} \geq A\sqrt n \, \|P\|_{[-1,1]}\,, \qquad A = \frac{2}{3\sqrt{210e}} = 0.0279\ldots \,,$$
for all polynomials of degree $n$ having all their zeros in the closed upper half-disk $D^+$. 

For integers $0 \leq k \leq n$ let ${\Cal F}_{n,k}^c$ be the set of all polynomials $P \in {\Cal P}_n^c$ having 
at least $n-k$ zeros in $D^+$. In this paper we prove an essentially sharp reverse Markov-type inequality for the 
classes ${\Cal F}_{n,k}^c$ extending the earlier results of Tur\'an and Komarov from the case $k=0$ to the cases 
$0 \leq k \leq n$.  

\head 2. New Results \endhead

Our main result in this paper is an essentially sharp reverse Markov-type inequality
for the classes ${\Cal F}_{n,k}^c$ extending the earlier results of Tur\'an and Komarov 
from the case $k=0$ to the cases $0 \leq k \leq n$. The upper bound of Theorem 2.1 below 
is quite a new result even in the case when the infimum is taken for polynomials $P \in {\Cal P}_n^c$ 
having at least $n-k$ zeros only in $[-1,1]$ rather than $D^+$. 

\proclaim{Theorem 2.1}
There are absolute constants $c_1 > 0$ and $c_2 > 0$ such that
$$c_1 \left(\frac{n}{k+1}\right)^{1/2} \leq 
\inf_{P}{\frac{\|P^{\prime}\|_{[-1,1]}}{\|P\|_{[-1,1]}}} \leq c_2 \left(\frac{n}{k+1}\right)^{1/2}$$
for all integers $0 \leq k \leq n$, where the infimum is taken for all $0 \not\equiv P \in {\Cal F}_{n,k}^c$ 
having at least one zero in $[-1,1]$. 
\endproclaim 

Theorem 2.1 follows from the results below. In fact, Corollary 2.3 below offers the explicit constant $c_1 = 1/808$ in 
a slightly modified form of the lower bound in Theorem 2.1. 

\proclaim{Theorem 2.2}
Let $1 \leq k \leq n/163000$. We have
$$\|P^\prime\|_{[-1,1]} \geq \frac{1}{202} \left( \frac{n-k}{8k} \right)^{1/2} \|P\|_{[-1,1]}$$
for all $P \in {\Cal F}_{n,k}^c$.
\endproclaim

\proclaim{Corollary 2.3}
Let $1 \leq k \leq n$. We have
$$\|P^\prime\|_{[-1,1]} \geq \max \left\{\frac 12,\frac{1}{808} \left( \frac{n-k}{k} \right)^{1/2} \right\} \|P\|_{[-1,1]}$$
for all $P \in {\Cal F}_{n,k}^c$ with at least one zero in $[-1,1]$.
\endproclaim

\proclaim{Theorem 2.4}
There is an absolute constant $c_1 > 0$ and there are polynomials 
$0 \not\equiv P = P_{n,k} \in {\Cal F}_{2n,2k}$ of the form 
$$P(x) = (x^2-1)^{n-k}R(x)\,, \qquad R \in {\Cal P}_{2k}\,,$$
such that 
$$\frac{\|P^\prime\|_{[-1,1]}}{\|P\|_{[-1,1]}} \leq c_2 \left( \frac nk \right)^{1/2}$$
for every $1 \leq k \leq n/2$.
\endproclaim

We remark that the upper bound of Theorem 2.1 does not remain valid if we replace the closed upper half-disk $D^+$ 
with the closed unit disk $D$ in the definition of ${\Cal F}_{n,k}^c$, not even in the case that $k=0$. 
Given $\varepsilon > 0$, let $m$ be the even integer for which $1/\varepsilon < m \leq 1/\varepsilon + 2$. 
We claim that for every $\varepsilon > 0$ and for every integer $n \geq 1$ there is a $P_n \in {\Cal P}_{mn}^c$ 
of degree $mn$ having all its zeros on the unit circle ${\partial D}$ such that 
$$\|P_n^\prime\|_{[-1,1]} \leq (1/\varepsilon + 2)^{1-\varepsilon}(mn)^{\varepsilon}\|P_n\|_{[-1,1]}\,.$$ 
To see this let $P_n \in {\Cal P}_{mn}^c$ be defined by  $P_n(z) := (z^m-1)^n$. 
Observe that $\|P_n\|_{[-1,1]} = 1$ (as $m$ is even), and the function 
$$|P_n^\prime(x)| = mn(1-x^m)^{n-1}|x|^{m-1}$$
achieves its maximum on $[-1,1]$ at the point $a \in (0,1)$, where i
$$a^m = \frac{m-1}{mn-1} \leq \frac 1n\,.$$ 
Hence
$$|P_n^\prime(a)| \leq mn a^{m-1} \leq mn n^{1/m-1} \leq mn^{\varepsilon} \leq m^{1-\varepsilon}(mn)^{\varepsilon}
\leq (1/\varepsilon + 2)^{1-\varepsilon}(mn)^{\varepsilon}\,.$$

\head 3. Lemmas \endhead

Our proof of Theorem 2.2 is based on the following two non-trivial results. Lemma 3.1
is proved in [GL-78] while the proof of Lemma 3.2 may be found in Section 7.2 of [BE-95].  

\proclaim{Lemma 3.1}
If $Q \in {\Cal F}_{n,0}^c$ and 
$$E_{\delta} := \left\{ x \in [-1,1]: \left| \frac{Q^\prime(x)}{Q(x)} \right| \leq n\delta \right\}\,, \qquad \delta > 0\,,$$
then 
$$m(E_{\delta}) < A\delta\,, \qquad \delta > 0\,,$$ 
where $A := 70e$ is a suitable choice.
\endproclaim

\proclaim{Lemma 3.2}
If $R \in {\Cal P}_k^c$ and
$$F_{\alpha} := \left\{ x \in {\Bbb R}: \left| \frac{R^\prime(x)}{R(x)} \right| \geq \alpha \right\}\,, \qquad \alpha > 0\,,$$
then
$$m(F_{\alpha}) \leq \frac{Bk}{\alpha}\,, \qquad \alpha > 0\,,$$
where $B := 8\sqrt{2}$ is a suitable choice.
\endproclaim

Our third lemma is a simple consequence of the Mean Value Theorem. 

\proclaim{Lemma 3.3}
If a function $P$ is differentiable on $[-1,1]$, 
$$\|P^\prime\|_{[-1,1]} \leq M \|P\|_{[-1,1]}\,,$$
and $x_0 \in [-1,1]$ is such that 
$$|P(x_0)| := \|P\|_{[-1,1]}\,, \tag 3.1$$
then
$$|P(y)| \geq \frac 12 \, \|P\|_{[-1,1]}\,, \qquad y \in [x_0-(2M)^{-1},x_0+(2M)^{-1}] \cap [-1,1]\,.$$
\endproclaim

To prove Theorem 2.4 we need the following two lemmas. Lemma 3.4 below is stated and proved as 
Theorem 2.1 in [E-19] by using deep results from [B-85] and [BE-94].  Recall that 
${\Cal P}_{n-k,k}$, $0 \leq k \leq n$, denotes the set of all algebraic polynomials 
with real coefficients, of degree at most $n$ having at least $n-k+1$ zeros at $0$.
 
\proclaim{Lemma 3.4} There are absolute constants $c_3 > 0$ and $c_4 > 0$ such that
$$c_3 \frac{n-k}{k} \leq \min_{0 \not \equiv P \in {\Cal P}_{n-k,k}}{\frac{\|P^{\prime}\|_{[0,1]}}{V_0^1(P)}}  
\leq \min_{0 \not\equiv P \in {\Cal P}_{n-k,k}}{\frac{\|P^{\prime}\|_{[0,1]}}{|P(1)|}} \leq c_4 \frac nk$$
for all integers $1 \leq k \leq n-1$. Here $c_3 = 1/12$ is a suitable choice.
\endproclaim

Lemma 3.5 below is stated and proved as Lemma 3.2 in [E-19].

\proclaim {Lemma 3.5}
Let $1 \leq k \leq n-1$ and let $S(x) := x^{n-k}R(x)$ with $R \in {\Cal P}_k$.
We have
$$|S(x)| \leq x^{(n-k)/2} \|S\|_{[0,1]}\,, \qquad x \in \left[0,1- \frac{10k}{n-k}\right]\,.$$
\endproclaim

Lemma 3.6 below follows simply from Lemma 3.5.

\proclaim {Lemma 3.6}
Let $1 \leq k \leq n/2$ and let $W(x) := (1-x)^{n-k}V(x)$ with $V \in {\Cal P}_k$.
We have
$$|y^{1/2}W(y)| < \|u^{1/2}W(u)\|_{[0,1]}\,, \qquad y \in \left[\frac{10(2k+1)}{n},1\right]\,.$$
\endproclaim

\demo{Proof of Lemma 3.6}
Replacing $n$ by $2n+1$ and $k$ by $2k+1$ in Lemma 3.6 we get obtain that 
$$|S(x)| \leq x^{n-k} \|S\|_{[0,1]}\,, \qquad x \in \left[0,1- \frac{10(2k+1)}{n}\right] \subset 
\left[0,1-\frac{10(2k+1)}{2n-2k}\right]\,, \tag 3.1$$
whenever $1 \leq k \leq n/2$ and $S(x) := x^{2n-2k}R(x)$ with $R \in {\Cal P}_{2k+1}$. 
Replacing the variable $x$ by $1-x$ in (3.1) yields that 
$$|S(x)| \leq (1-x)^{n-k} \|S\|_{[0,1]}\,, \qquad x \in \left[\frac{10(2k+1)}{n},1\right]\,, \tag 3.2$$
whenever $1 \leq k \leq n/2$ and $S(x) := (1-x)^{2n-2k}R(x)$ with $R \in {\Cal P}_{2k+1}$.
Now let $1 \leq k \leq n/2$ and let $W(x) := (1-x)^{n-k}V(x)$ with $V \in {\Cal P}_k$. 
Applying (3.2) to $S$ defined by 
$$S(x) = xW(x)^2 = (1-x)^{2n-2k}(xV(x)^2)\,, \qquad V \in {\Cal P}_k\,,$$ 
we get the conclusion of the lemma. 
\qed \enddemo

\head 4. Proof of the Theorems \endhead

\demo{Proof of the Theorem 2.2}
Let $0 \not \equiv P \in {\Cal F}_{n,k}^c$, that is, $P = QR$, where $P \in {\Cal F}_{n-k,0}^c$ and 
$R \in {\Cal P}_k^c$. We have
$$\frac{P^\prime}{P} = \frac{Q^\prime}{Q} + \frac{R^\prime}{R}\,. \tag 4.1$$
By Lemma 3.1 we have
$$m(E_{\delta}) < A\delta, \qquad \delta > 0, \qquad A := 70e\,, \tag 4.2$$
where 
$$E_{\delta} := \left\{ x \in [-1,1]: \left| \frac{Q^\prime(x)}{Q(x)} \right| \leq (n-k)\delta \right\}\,, 
\qquad \delta > 0\,. \tag 4.3$$
By Lemma 3.2 we have
$$m(F_{\delta}) \leq B\delta, \qquad \delta > 0, \quad B := 8\sqrt{2}\,, \tag 4.4$$
where
$$F_{\delta} := \left\{ x \in [-1,1]: \left| \frac{R^\prime(x)}{R(x)} \right| \geq \frac{k}{\delta} \right\}\,, 
\qquad \delta > 0\,. \tag 4.5$$
Now we choose $\delta > 0$ such that
$$\frac{k}{\delta} = \frac 12 \, (n-k)\delta\,, \tag 4.6$$
that is, 
$$\delta := \left( \frac{2k}{n-k} \right)^{1/2}\,. \tag 4.7$$
Then, combining (4.1)--(4.7), we can deduce that
$$\split \left| \frac{P^\prime(x)}{P(x)} \right| \geq &
\left| \frac{Q^\prime(x)}{Q(x)} \right| - \left| \frac{R^\prime(x)}{R(x)} \right| \cr 
\geq & (n-k)\delta - \frac{k}{\delta} \geq \frac 12 \, (n-k)\delta 
\geq \frac 12 \, (n-k)\left( \frac{2k}{n-k} \right)^{1/2} \cr 
\geq & \left( \frac 12 \, (n-k)k \right)^{1/2}\,, \qquad x \in [-1,1] \setminus H_{\delta}\,, \cr \endsplit \tag 4.8$$
where $H_{\delta} := E_{\delta} \cup F_{\delta}$
with 
$$m(H_{\delta}) \leq (A+B)\delta\,. \tag 4.9$$
Note that 
$$1 \leq k \leq \frac{n}{163000}$$
implies that
$$(A+B)\delta = (A+B)\left( \frac{2k}{n-k} \right)^{1/2} 
\leq (70e + 8\sqrt{2})\left( \frac{4k}{n} \right)^{1/2} < 1\,. \tag 4.10$$
Assume now to the contrary of the theorem that  
$$\|P^\prime\|_{[-1,1]} \leq (2(A+B)\delta)^{-1} \|P\|_{[-1,1]}\,. \tag 4.11$$
Choose $x_0 \in [-1,1]$ such that (3.1) holds. Observe that (4.9), (4.10), and (4.11) imply that there is a 
$$y \in [x_0-(A+B)\delta,x_0+(A+B)\delta] \cap [-1,1] \tag 4.12$$
such that 
$$y \in [-1,1] \setminus H_{\delta}\,. \tag 4.13$$ 
Using Lemma 3.3 with 
$M := (2(A+B)\delta)^{-1}$ and recalling (4.12) we obtain 
$$|P(y)| \geq \frac 12 \, \|P\|_{[-1,1]}\,. \tag 4.14$$
Combining (4.13), (4.8) and (4.14) we obtain
$$\split \|P^\prime\|_{[-1,1]} \geq & |P^\prime(y)| > \left( \frac 12 \, (n-k)k \right)^{1/2}|P(y)| \cr 
\geq & \left( \frac 12 \, (n-k)k \right)^{1/2} \frac 12 \, \|P\|_{[-1,1]} \cr
> & (2(A+B)\delta)^{-1} \|P\|_{[-1,1]} \,, \cr \endsplit$$
which contradicts (4.11). Hence (4.11) is impossible and the proof is finished. 
\qed \enddemo

\demo{Proof of Corollary 2.3}
Let $1 \leq k \leq n$. Suppose $0 \not\equiv P \in {\Cal F}_{n,k}^c$ has at least one zero in $[-1,1]$. Choose 
$a,b \in [-1,1]$ such that $P(a) = 0$, and $|P(b)| = \|P\|_{[-1,1]}$. By the Mean Value Theorem 
there is a $c \in (-1,1)$ between $a$ and $b$ such that
$$\|P^\prime\|_{[-1,1]} \geq |P^\prime(c)| = \left| \frac{P(b)-P(a)}{b-a} \right| \geq \frac 12 \|P\|_{[-1,1]} \tag 4.15$$
If $\displaystyle{1 \leq k \leq \frac{n}{163000}}$, the result follows from Theorem 2.2 and (4.15). 
If $\displaystyle{\frac{n}{163000} < k \leq n}$, then
$$\frac{1}{808} \left( \frac{n-k}{k} \right)^{1/2} \leq \frac 12\,,$$ 
and the result follows simply from (4.15). 
\qed \enddemo

\demo{Proof of Theorem 2.4}
Let $1 \leq k \leq n/2$. By (the upper bound of) Lemma 3.4 there is an absolute constant 
$c_4 > 0$ and there are polynomials 
$0 \not\equiv Q = Q_{n,k} \in {\Cal P}_{n-k,k}$ such that 
$$\frac{\|Q^\prime\|_{[0,1]}}{\|Q\|_{[0,1]}} \leq c_4 \frac nk\,. \tag 4.16$$
Let $0 \not\equiv R(x) = R_{n,k}(x) = Q(1-x)$. Obviously 
$R$ is of the form 
$$R(x) = (1-x)^{n-k+1}U(x)\,, \qquad U \in {\Cal P}_{k-1}^c\,,$$
$R^\prime$ is of the form 
$$R^\prime(x) = (1-x)^{n-k}V(x)\,, \qquad V \in {\Cal P}_{k-1}^c\,, \tag 4.17$$
Let $0 \not\equiv P = P_{n,k}$ be defined by $P(x) := R(x^2)$. Observe that $P$ is of the form 
$$P(x) = (1-x^2)^{n-k+1}U(x)\,, \qquad U \in {\Cal P}_{2k-2}^c\,,$$
hence $P \in {\Cal F}_{2n,2k}$.
It is also clear that 
$$\|P\|_{[-1,1]} = \|R\|_{[0,1]} = \|Q\|_{[0,1]}  \tag 4.18$$
and 
$$P^\prime(x) = 2xR^\prime(x^2)\,. \tag 4.19$$ 
Let $y := x^2$. Using (4.19), (4.17), and Lemma 3.6, we obtain  
$$|P^\prime(x)| = |xR^\prime(x^2)| = |2y^{1/2}R^\prime(y)| < \|2u^{1/2}R^\prime(u)\|_{[0,1]} = \|P^\prime\|_{[-1,1]}$$
or every $y = x^2 \in [10(2k+1)/n,1]$, and hence there is an 
$$a \in \left[0, \left( \frac{10(2k+1)}{n} \right)^{1/2} \right] \tag 4.20$$
such that 
$$|P^\prime(a)| = \|P^\prime\|_{[0,1]}\,.  \tag 4.21$$
Combining (4.21), (4.19), (4.20), (4.16), and (4.18) we obtain 
$$\split \|P^\prime\|_{[-1,1]} & = \|P^\prime\|_{[0,1]} = |P^\prime(a)| = |2aR^\prime(a^2)| \cr
& \leq 2\left( \frac{10(2k+1)}{n} \right)^{1/2}\|R^\prime\|_{[0,1]} = 
2\left(\frac{10(2k+1)}{n}\right)^{1/2}\|Q^\prime\|_{[0,1]} \cr 
& \leq 2\left(\frac{10(2k+1)}{n}\right)^{1/2}c_4 \frac nk \|Q\|_{[0,1]} \cr 
& \leq c_2\left(\frac{n}{k}\right)^{1/2} \|Q\|_{[0,1]} = c_2\left(\frac{n}{k}\right)^{1/2} \|P\|_{[-1,1]}\,. \cr \endsplit$$
\qed \enddemo

\demo{Proof of Theorem 2.1}
The case that $k=0$ is the result of Komarov [K-19] mentioned in the Introduction, so we may assume that 
$1\leq k \leq n$. in which cases the lower bound of the theorem follows immediately from Corollary 2.3. 
Let $f(n,k)$ defined by
$$f(n,k) := \min_{0 \not \equiv P \in {\Cal F}_{n,k}^c}{\frac{\|P^{\prime}\|_{[-1,1]}}{\|P\|_{[-1,1]}}}\,.$$
Observe that for a fixed positive integer $n$ the function $f(n,k)$ is decreasing on the set of integers 
$0 \leq k \leq n$, and for a fixed nonnegative integer $k$ the function $f(n,k)$ is decreasing on the set of  
integers $n \geq k$. So it is sufficient to show the upper bound of the theorem only for even numbers 
$n = 2\nu$ and $k = 2\kappa$ satisfying $1 \leq \kappa \leq \nu/2$ in which cases the upper bound of the 
theorem follows from Theorem 2.4.
\qed \enddemo

\head 5. Acknowledgement \endhead
The author thanks Szil\'ard R\'ev\'esz for checking the details of the proofs
in a first draft of this paper and for helpful discussions.

\enddocument